# Solution Bounds, Stability and Attractors' Estimates of Some Nonlinear Time-Varying Systems

Mark A. Pinsky and Steve Koblik

**Abstract.** Estimation of solution norms and stability for time-dependent nonlinear systems is ubiquitous in numerous applied and control problems. Yet, practically valuable results are rare in this area. This paper develops a novel approach, which bounds the solution norms, derives the corresponding stability criteria, and estimates the trapping/stability regions for a broad class of the corresponding systems. Our inferences rest on deriving a scalar differential inequality for the norms of solutions to the initial systems. Utility of the Lipschitz inequality linearizes the associated auxiliary differential equation and yields both the upper bounds for the norms of solutions and the relevant stability criteria. To refine these inferences, we introduce a nonlinear extension of the Lipschitz inequality, which improves the developed bounds and estimates the stability basins and trapping regions for the corresponding systems. Finally, we conform the theoretical results in representative simulations.

## 1. INTRODUCTION.

We are going to study a system defined by the equation

(1) $\qquad \dot{x} = A(t)x + f(t,x) + F(t), \ t \in [t_0, \infty), \ x \in \mathbb{R}^n, \ f(t,0) = 0$

where matrix $A \in \mathbb{R}^{n \times n}$ and functions $f : [t, \infty) \times \mathbb{R}^n \to \mathbb{R}^n$ and $F : \mathbb{R} \to \mathbb{R}^n$ are continuous and bounded, and $A(t)$ is also continuously differentiable. It is assumed that the solution, $x(t, x_0)$, to the initial value problem, $x(t_0, x_0) = x_0$, is uniquely defined for $\forall t \geq t_0$. Note that the pertained conditions can be found, e.g., in [1] and [2].

We also examine the solutions to a homogeneous counterpart to (1)
(2) $\qquad \dot{x} = A(t)x + f(t,x)$

Development of efficient stability criteria for $x = 0$ solution to (2) is essential in numerous applied problems. There are two main approaches to this problem: the Lyapunov functions method, see for instance [2] and [3], and the first approximation methodology, see e.g. recent reviews [4] and [5] as well as [14] and [15] for additional references and historical perspectives. The former approach, for instance, is widespread in control literature, see [2] and [6]- [13] and additional references therein. However, adequate Lyapunov functions are rare, especially for time dependent and nonlinear systems. The latter approach delivers sufficient stability criteria under two conditions, see [14], and [15]. The first is the Lipschitz condition
(3) $\qquad \|f(x,t)\| \leq l(t)\|x\|, \ \forall x \in \Omega, \ \forall t \geq t_0$

where $\Omega$ is a neighborhood of $x = 0$ and $l(t) \leq \hat{l}$ be continuous function. The second condition,
(4) $\qquad \|W(t,t_0)\| \leq N e^{-\lambda(t-t_0)}, \ \forall t \geq t_0, \ 0 < N, \lambda = const$

bounds the growth rate of the transition matrix $W(t,t_0) = W(t)W^{-1}(t_0)$, where $W(t)$ is the fundamental solution matrix of

M.A. Pinsky is with the Department of Mathematics & Statistics, University of Nevada. Reno, Reno NV 89521 USA (e-mail: pinsky@unr.edu). Steve Koblik is independent consultant in scientific computing (e-mail: Stevekoblik8110@comcast.net)

.

(5) $$\dot{x} = A(t)x$$

Inequality (4) comprises necessary and sufficient conditions for asymptotic/exponential stability of (5), e.g. [3] and [14]. Consequently, it was shown that $x = 0$ solution to (2) is exponentially stable if (3), (4), and the following condition,

(6) $$N\bar{l} - \lambda < 0$$

are satisfied [14], [15]. A more flexible sufficient condition,

(7) $$\|W(t,t_0)\| \leq N \exp \int_{t_0}^{t} p(s)ds, \ \forall t > t_0 > 0$$

which reduces to (4) for $p(t) = const$ was introduced in [16], see also [4] and [5]. In turn, (3) and (7) provide asymptotic stability of (2) if [4],

(8) $$\limsup_{t \to \infty} (1/t) \int_{t_0}^{t} p(s)ds + N\bar{l} < 0$$

While the existence of (4) is acknowledged under some broad conditions [14], to our knowledge, there were no attempts to apply bound (4) or (7) to stability analysis of the systems of practical importance. Furthermore, it was shown, e.g. in [17], that the time-histories of different estimates of the bounds of the Euclidian norms for the second order fundamental matrix, i.e. $\|W(t)\| = \|\exp At\|$, $A = const$, can drastically diverge from each other and the exact values of $\|\exp At\|$. This raises concern of the practical value of the listed above sufficient stability criteria.

An attempt to escape the utility of prior bounds on $\|W(t,t_0)\|$ in stability analysis of (2) was undertaken in [18]. However, authentication of the developed stability conditions for relatively complex systems can present a challenging task for this approach as well.

The problem of estimating the norms of solutions to (1) subject to (3) and (4) was reviewed in [14] and [15]. This problem has some connection to input-to-state stability, which was mainly studied in the context of the Lyapunov function methodology in control literature, see [3] and [19] for further details and references.

This paper derives a novel scalar differential inequality for the norms of solutions to (1) or (2), which collapses the dimension of the original estimation problem to one. Due to the comparison principle [3], this inequality further reduces to the analysis of solutions to the auxiliary first order scalar nonlinear equation with variable coefficients. Utility of prior bounds on $\|W(t,t_0)\|$ is naturally voided in our study, which splits in two parts. The first adopts (3) which linearizes the auxiliary equation and yields the solution bounds and corresponding stability criteria. The second part introduces a nonlinear extension of the Lipschitz inequality. This sharpens the bounds of solutions and estimates the trapping/instability regions for nonlinear equations with time dependent coefficients. Our inferences are validated in simulations of the Van der- Pol like model, which includes a time dependent linear block and oscillatory external force.

This paper is organized as follows. The next section derives the pivotal differential inequality, the subsequent section develops solution bounds and some stability criteria via utility of the Lipschitz inequality, section 4 introduces a nonlinear extension of the Lipschitz inequality and develops its applications, section 5 includes the simulation results, and section 6 concludes this study.

2. Differential Inequality for Solution Norms

The application of variation of parameters lets us write (1) as [3]

$$x(t,x_0) = W(t)W^{-1}(t_0)x_0 + W(t)\int_{t_0}^{t} W^{-1}(\tau)\big(f(\tau,x(\tau)) + F(\tau)\big)d\tau$$

where $W(t)$ is frequently normalized to satisfy the condition, $W(t_0) = I$, $I$ is the identity matrix. In section 5 we present normalization of $W(t)$, which is more natural for our studies and, hence, used in our simulations. Presently, we only assume that $\|W(t_0)\| = 1$. The last equation yields the following inequality,

(9) $$\|x(t,x_0)\| \leq \|W(t)\|\|W^{-1}(t_0)x_0\| + \|W(t)\|\int_{t_0}^{t}\|W^{-1}(\tau)\|\|f(\tau,x(\tau)) + F(\tau)\|d\tau$$

Now we attempt to match (9) with a first order differential inequality and the associate initial condition as follows,

(10) $$D^+\|x\| \leq p(t)\|x\| + k(t)\|f(t,x) + F(t)\|$$
$$\|x(t)\|_{t=t_o} = X_0$$

where $D^+$ is Dini's upper right-hand derivative in $t$ [3]. In turn, the application of variation of parameters to (10) yields,

$$\|x\| \leq e^{\zeta}\left(X_0 + \int_{t_0}^{t} e^{-\zeta} k(\tau)\psi(\tau)d\tau\right)$$

where $\zeta = \int_{t_0}^{t} p(s)ds$ and $\psi = \|f(\tau,x(\tau)) + F(\tau)\|$. Comparison of (9) with the last formula returns,

$$\|W(t)\| = \exp\left(\int_{t_0}^{t} p(s)ds\right), \quad X_0 = \|W^{-1}(t_0)x_0\|$$

Then, we can write that

$$\int_{t_0}^{t} e^{-\zeta} k(\tau)\psi(\tau)d\tau = \int_{t_0}^{t}\left(\|W(\tau)\|\|W^{-1}(\tau)\|\right)\psi(\tau)/\|W(\tau)\|d\tau$$

Hence, $k(t) = \|W(t)\|\|W^{-1}(t)\| = \sigma_{max}(W)/\sigma_{min}(W)$ is the running condition number of $W(t)$, and $\sigma_{max}$ and $\sigma_{min}$ are maximal and minimal singular values of $W(t)$. We assume that, $\forall t \geq t_0$, $\sup \sigma_{max} \leq \bar{\sigma} = const$ and $\inf \sigma_{min} \geq \underline{\sigma} = const$ and $\sigma_{max}$ is unique. This, in particular, assures that $k(t)$ is continuous and $|k(t)| \leq \hat{k} = const$. We also write that

(11) $$p(t) = d\left(\ln\|W(t)\|\right)/dt$$

Hence, $p(t)$ is continuous since we assumed that $\sigma_{max}$ is unique.

To write (10) in the standard form, we introduce a nonlinear extension of Lipschitz inequality,

(12) $$\|f(x,t)\| \leq L(t,\|x\|), \quad \forall x \in \Omega, \ \forall t \geq t_0$$

where $\Omega$ is a neighborhood of $x = 0$ and $L$ is a function continuous in $t$ and $\|x\|$. In section 4 we review how to devise (12) for some broad sets of functions. Clearly, (12) reduces to (3) if $L$ is linear in $\|x\|$. Finally, we write (10) in the following form

(13) $$D^+\|x\| \leq p(t)\|x\| + k(t)\left(L(t,\|x\|) + \|F(t)\|\right)$$
$$\|x(t_0)\| = \|W^{-1}(t_0)x_0\|$$

Interestingly, the sensitivity of (5) to perturbations is scaled by $k(t)$ and, hence, varies in time. Due to the comparison principle [3], solutions of (13) are bounded by the corresponding solutions to the associated first order auxiliary differential equation, which is analyzed in section 4. In the following section we linearize (13) through utility of (3). This lets us establish some upper bounds on $\|x(t,x_0)\|$ and derive the corresponding stability criteria.

3. Linearized Auxiliary Equation

Using (3), we write (13) as a linear inequality,

$$D^+ \|x\| \leq (p(t) + k(t)l(t))\|x\| + k(t)\|F(t)\|, \quad \forall x \in \Omega, \ \forall t \geq t_0 \tag{14}$$

$$\|x(t_0)\| = \|W^{-1}(t_0)x_0\|$$

Due to the comparison principle [3], solutions of (14) are bounded by the consistent solutions of the corresponding auxiliary differential equation. We mention that our previous assumptions assure that the initial value problem for this equation possesses a unique solution. This leads to the following,

**Theorem 1**. Assume that (3) is intact and, due to our previous assumptions, $k$ and $\|F\|$ are continuous and bounded, and $p$ is continuous. Then,

$$\|x(t, x_0)\| \leq \|x_h(t, x_0)\| + \|x_{nh}(t, x_0)\|, \quad \forall x \in \Omega, \ t_0 \leq t < \infty \tag{15}$$

where

$$\|x_h\| = \|W(t)\| \|W^{-1}(t_0)x_0\| \exp\left(\int_{t_0}^{t} k(s)l(s)ds\right) \tag{16}$$

and

$$\|x_{nh}\| = \int_{t_0}^{t} \theta(t, \tau)k(\tau)\|F(\tau)\| d\tau \tag{17}$$

where transition function, $\theta(t, \tau) = \exp\left(\int_{\tau}^{t} (p(s) + k(s)l(s))ds\right)$.

Clearly, due to the assumptions made in previous section, the integrals in (16) and (17) exists for $t_0 \leq t < \infty$, which assures this statement.

Hence, the problems of assessing the asymptotic/exponential stability of $x = 0$ solution to (2) or boundedness of solutions to (1) with $x_0 \in \Omega$ are simplified and comprised in evaluation of the matching problems to the auxiliary linear first order homogeneous/nonhomogeneous equations, and assuring that $x(t, x_0) \in \Omega, \forall t \geq t_0$. The latter, in turn, can be inferenced if either $\|x_0\|$ or both $\|x_0\|$ and $\|F(t)\|$ are sufficiently small. The necessary and sufficient conditions for various types of stability of linear equations are known, see, e.g. [14] and are recently reviewed, e.g., in [20], where additional references can be found. Application of these conditions to our first order linear auxiliary equation yields the matching stability criteria for $x = 0$ solution to nonlinear equation (2). Below we formulate some of most explicit stability/ boundedness conditions.

**Corollary 1.** Let $\|W(t)\|k(t)l(t) + d\|W(t)\|/dt < 0$ for $\forall t > t_1 \geq t_0$. Then $x = 0$ solution to (2) is asymptotically stable and $\|x_h\| \to 0$ monotonically for $t > t_1$.

In fact, due to (16), $\|x_0\|$ can be chosen such that $\|x(t, x_0)\| \in \Omega, \forall t \in [t_1, t_0]$. Then, the statement is assured since the derivative of the right-side of (16) is negative for $t > t_1$.

**Corollary 2.** Let $\sup(p(t) + k(t)l(t)) < -v_1, \ v_1 > 0, \ \forall t > t_1 \geq t_0$. Then, $x = 0$ solution of (2) is exponentially stable. If, in addition, $\|F(t)\| \leq F_0$ and both $F_0$ and $\|x_0\|$ are sufficiently small, then the corresponding solutions to (1) are bounded.

This statement directly follows from (15) - (16).

To formulate less conservative stability criteria, we evoke the definitions of characteristic and Lyapunov exponents, which determine the fate of solutions to (2) and (5) for $t \to \infty$, see e.g., [4]. The characteristic exponents, $\chi_i(x(t, x_0)) = \limsup_{t \to \infty} t^{-1} \ln \|x(t, x_0)\|, \ i \leq n$, assess the rate of exponential growth/decay of $\|x(t, x_0)\|$ if $t \to \infty$. For a linear system, the Lyapunov exponents are defined as, $\bar{\chi}_i = \limsup_{t \to \infty} t^{-1} \ln \sigma_i(\|W(t)\|), \ i \leq n$, where $\sigma_i$ are the singular values of $W(t)$. For linear systems, maximal characteristic and Lyapunov exponents are matched, see e.g. [4].

Firstly, we notice that

$$\limsup_{t\to\infty} t^{-1}\int_{t_0}^{t} p(s)ds = \limsup_{t\to\infty} t^{-1}\ln\|W(t)\| = \limsup_{t\to\infty} t^{-1}\ln\sigma_{\max}(t) = \bar{\chi}_{\max}$$

where $\bar{\chi}_{\max}$ is the maximal Lyapunov exponent for (5) and $\sigma_{\max}$ is the maximal singular value of $W(t)$. Hence, (5) is exponentially stable if $\bar{\chi}_{\max} < 0$ and uniformly exponentially stable if $\inf \sigma_{\min} \geq \underline{\sigma} = const$, as we assumed in the previous section. Next, we can calculate the upper bound for characteristic exponents for (2), i.e. $\hat{\chi} \geq \chi_i$, using (16),

$$\hat{\chi} = \limsup_{t\to\infty} t^{-1}\left(\int_{t_0}^{t} k(s)l(s)ds + \ln\|W(t)\|\right)$$

Then, it follows that $|\hat{\chi}| < \infty$ since $\|W(t)\| \leq \exp\int_{t_0}^{t}\|A(s)\|ds$, see e.g. [14] and $\|A\|, |k|, |l| < \infty$ by the prior assumptions. Now, let us set $\hat{\chi} = \bar{\chi}_{\max} + \chi_*$, where $\chi_* = \limsup_{t\to\infty} t^{-1}\int_{t_0}^{t} k(s)l(s)ds$, then we convey,

**Corollary 3**. Solution $x = 0$ to (2) is exponentially stable if $\hat{\chi} < 0$. If, in addition, $\inf \sigma_{\min} \geq \underline{\sigma} = const$, then this solution is uniformly exponentially stable. Indeed, it follows from (16) that

$$\|x(t,x_0)\| \leq \|W^{-1}(t_0)\|\|x_0\|\exp\int_{t_0}^{t}(p(s)+k(s)l(s))ds, \ \forall t \geq t_0$$

The last multiple in the above formula decays exponentially fast since $\hat{\chi} < 0$, and note also that $\|W^{-1}(t_0)\| = \sigma_{\min}^{-1}(t_0) \leq \underline{\sigma}^{-1}$.

**Corollary 4**. Let $\hat{\chi} < 0$, $\inf \sigma_{\min} \geq \underline{\sigma} = const$, $\|F(t)\| \leq F_0$, and $\|x_0\|$ and $F_0$ be sufficiently small. Then the corresponding solutions to (1) are bounded in norm $\forall t \geq t_0$.

In fact, in this case the transition function for our auxiliary equation, $\theta(t,\tau) \leq Me^{-\lambda(t-\tau)}$, where $M > 0$, $-\lambda = \hat{\chi} + \varepsilon$ and $0 < \varepsilon$ can be chosen sufficiently small, see e.g., [21], pp.100-101. Then, due to variation of parameters, we get,

$$\|x_{nh}(t,x_0)\| \leq F_0 M \int_{t_0}^{t} e^{-\lambda(t-\tau)}d\tau = (F_0 M / \lambda)(1-\exp(-\lambda(t-t_0)))$$

Next, due to Corollary 3, $\|x(t,x_0)\| \in \Omega$ if $\|x_0\|$ and $F_0$ are sufficiently small.

*Remark 1*. We notice that the application of stability criteria (developed for a scalar linear system in [20]) to our auxiliary equation (14) leads to somewhat less conservative stability criteria of nonlinear equation (2). This lets us to replace estimation of the characteristic exponents by computing the pertained moving averages over sufficiently long-time intervals. These types of augments of the above statements are left out of this paper. Yet, in our view, numerical estimation of the characteristic exponents and computation of the corresponding time-averages seem to be quite similar in effort.

### 4. Nonlinear Lipschitz Inequality and its Applications
### 4.1 Nonlinear Lipschitz Inequality

Though Lipschitz condition was used extensively in stability theory, e.g., [3]- [5], [14], and [15], its applications frequently lead to over conservative estimations, which also evoke dependence of the Lipschitz constant upon the size of the pertaining neighborhood, i.e., $l = l(\Omega)$. A rigorous assessment of the last relation can present a complex task, which is frequently avoided. Yet, this can affect the practical value of the pertained results. Additionally, utility of (3) linearizes (10) and abates representation of intrinsically nonlinear phenomena.

To temper these problems, we consider a nonlinear extension of Lipschitz inequality, i.e., (12). In principle, a somewhat conservative form of (12) can be readily derived in many cases. For instance, for polynomial vector fields or ones, which can be approximated by interpolating/Taylor's polynomials with globally bounded Lagrange error terms, (12) converts to a global inequality. In these cases, (12) can be attained in polynomial form, e.g., by successive applications to every addition of these approximations of the following inequalities: $\|f\|_2 \leq \|f\|_1$ and $|x_m^k| \leq \|x\|_2^k$, $m = 1,...,n$, where $x_m$ is the $m$-th component of $x$.

If the error term in the polynomial approximation of $f$ is bounded for $x \in \Omega$, then (12) is validated in the same neighborhood. Yet, often such nonlinear inequalities turns out to be less conservative than (3) in extended neighborhoods of $x = 0$ solution to (2) and lead to better apprehension of nonlinear nature of these equations.

### 4.2 Solution Bounds and Attractors' Estimates

Utility of (13) frequently sharpens the corresponding estimates and lessens their dependence upon the size of the pertaining neighborhood $\Omega$ but leads to analysis of the IVP for a nonintegrable and nonlinear first order ode with variable coefficients,

(18)
$$D^+ X = p(t) X + k(t) L(t, X) + k(t) \|F(t)\|$$
$$X(t_0) = X_0$$

where $X = \|x\|$ and $X_0 = \|W^{-1}(t_0) x_0\|$. Like above, we assume that $L(t, X)$ and $\|F(t)\|$ are continuous in their variables and (18) admits a unique solution, $\forall X_0 \geq 0$, $t \geq t_0$. Note that the last two terms in the right side of (18) are nonnegative, whereas $p(t) X$ can be either positive or negative or switch the sign for certain values of $t$.

Next, we reveal the correspondence between the solutions to scalar equation (18), i.e $X(t, X_0)$ and the norms of solutions to (1) or (2), i.e. $\|x(t, x_0)\|$. For this sake we define the close set of initial vectors to equation (1), i.e., $x_0 \in \omega(t_0, X_0)$ as follows,

$$\omega(t_0, X_0) \mid \|W^{-1}(t_0) x_0\| \leq X_0$$

These sets are bounded by concentric ellipsoids, $\omega_\gamma(t_0, X_0)$, which are centered at $x_0 = 0$. This leads to the following,

**Theorem 2**. Let $x_0 \in \omega(t_0, \hat{X})$, then,

(19)
$$\|x(t, x_0)\| \leq X(t, \hat{X}), \ \forall t \geq t_0$$

where $x(t, x_0)$ is a solution to (1).

In fact, $\|x(t, x_0)\| \leq X(t, X_0)$. Let $X_0 \leq \hat{X}$, then, $X(t, X_0) \leq X(t, \hat{X})$, since the solutions to (18) do not intersect due to uniqueness, which implies (19).

Inequality (19) facilitates estimation of the trapping/stability regions for (18), which, in turn, lets us to estimate the corresponding regions for systems (1) and (2), respectively.

**Corollary 5.** Suppose (18) assumes one of the following two conditions:
1. $F = 0$ and $\lim_{t \to \infty} X(t, \hat{X}) = 0$, then $x = 0$ is the asymptotically stable fixed solution to (2) and $\omega(t_0, \hat{X})$ is enclosed in its stability basin.
2. $F \neq 0$ and $X(t, \hat{X}) \leq d < \infty$, and $x_0 \in \omega(t_0, \hat{X})$, then $\|x(t, x_0)\| \leq X(t, \hat{X})$, where $x(t, x_0)$ is a solution to (1), i.e., $\omega(t_0, \hat{X})$ is enclosed in the trapping region of solutions to (1).

This corollary directly follows from (19).

Clearly, the best estimates of the trapping/stability regions correspond to the pertained value of $\max \hat{X}$, which obey the above statement. These values can be readily assessed in simulations of a scalar equation (18) especially since $X(t, X_0)$ is an increasing function in $X_0$ for every fixed $t$. Below, we formulate two complementary analytical approaches to estimation of the corresponding values of $\hat{X}$, which also let us better understand the structure of solutions to (18).

### 4.3 Explicit Estimation of Attractors

Let us turn (18) into a more conservative, but time invariant and integrable form

(20)
$$D^+ X = \hat{p}X + \hat{k}\hat{L}(X) + \hat{k}\hat{F} = Q(X, \hat{F})$$
$$X(0) = X_0$$

where we assumed that,

(21)
$$\sup k(t) = \hat{k}, \sup p(t) = \hat{p}, \sup_t L(t, X) = \hat{L}(X), \sup \|F(t)\| = \hat{F}, \forall t \geq t_0$$

and $\hat{k}, \hat{p}, \hat{F}, \hat{L}(X) < \infty$, $\hat{L}(0) = 0$. We also assume that (20) has a unique solution for the pertained initial values, which can be extended for $\forall t \geq t_0$.

As is known, the positive roots of equation

(22)
$$Q(X, \hat{F}) = 0$$

i.e., $X = d_i, i = 1, ..., I$, split the solutions to (20) into subsets with different behavior for $t \to \infty$. Subsequently, equation

(23)
$$\|W^{-1}(t_0)x_0\| = d_i$$

corresponds these point-wise boundaries to $n-1$-dimensional ellipsoids in the phase space of (1) or (2), i.e., $\omega_\gamma(t_0, d_i)$, which estimate the trapping/stability regions for the conforming systems. Below, we review the application of this procedure to some relatively simple, but characteristic, cases.

Assume, firstly, that $\hat{p} > 0$, then $Q(X, \hat{F}) \geq 0$. Hence, if $\hat{F} \neq 0$, then $\forall X_0 \geq 0$, $X(t, X_0) \to \infty$ monotonically if $t \to \infty$. If $\hat{F} = 0$, then $\forall X_0 > 0$, $X(t, X_0) \to \infty$ monotonically if $t \to \infty$. However, in these both case, the corresponding norms of solutions to (1) or (2) can either approach positive infinity, zero, or remain to be bounded.

Assume now that $F = 0$, $\hat{p} < 0$, and $\hat{L}(X)$ is monotonically increasing for $X > 0$. Then, (22) has no or one root $d > 0$, i.e., $Q(d, 0) = 0$. Note, that the former case, where $Q(X, 0) \geq 0$, was already conversed above. In the latter case, the fixed solution, $X = d$, can be either stable or unstable. This yields the following,

**Theorem 3**. Let $F = 0$, $\hat{p} < 0$, $\hat{L}(X)$ be a monotonically increasing function, and $X = d$ be a unique fixed solution to (20).

1. If this solution is unstable, then $x = 0$ solution to (2) is asymptotically stable. Furthermore, if the initial vector, $x_0 \in \omega(t_0, d)$, then $\|x(t, x_0)\| \leq d$, $\forall t \geq t_0$ and $\lim_{t \to \infty} \|x(t, x_0)\| = 0$, i.e., $\omega(t_0, d)$ belongs to the stability basin of $x = 0$ solution to (2).

2. If $X = d$ is a stable solution to (20), then $\|x(t, x_0)\| \leq d$, $\forall x_0 \in \omega(t_0, d)$ and $\limsup_{t \to \infty} \|x(t, x_0)\| \leq d$, $\forall x_0 \geq 0$, i.e., $\omega(t_0, d)$ belongs to the trapping region of $x = 0$ solution to (2).

The proof of this theorem immediately follows from (19) and the assessment of the behavior of solutions to (20) in these two cases. Note that the second case also includes occasions when $\hat{L}(X)$ and, hence, $f(t, x)$ are not Lipschitz continuous at $X = 0$ and $x = 0$ respectively. Such systems are excluded from this paper since the uniqueness of solutions at the relevant points for these equations are not warrantied.

Next, we assume that $\hat{F} \neq 0$, $\hat{p} < 0$, and $\hat{L}(X)$ is a monotonically increasing function. Then, (22) has no or two positive roots. The former case was already conversed above; whereas in the latter scenario, the roots $d_i, i = 1, 2$ can be either equal or distinct. This comprises,

**Theorem 4.** Let $\hat{F} \neq 0, \hat{p} < 0$, $\hat{L}(X)$ be a monotonic function, and $x(t, x_0)$ be a solution to (1).

1. Assume that $d_1 > d_2 > 0$ are the roots of (22) corresponding to unstable and stable fixed solutions of (20). Then:

$$\|x(t,x_0)\| \leq \begin{cases} d_1, & \forall x_0 \in \omega(t_0,d_1) \\ d_2, & \forall x_0 \in \omega(t_0,d_2) \end{cases}, \forall t \geq t_0,$$

$$\limsup_{t \to \infty} \|x(t,x_0)\| \leq d_2, \forall x_0 \in \omega(t_0,d_1)$$

2. Let $d_1 = d_2 = d$, then $\|x(t,x_0)\| \leq d, \forall x_0 \in \omega(t_0,d), t \geq t_0$.

The proof of this statement immediately follows from (19) and the analysis of the behavior of solutions to (20) in the relevant cases.

Obviously, (22) can admit more than two positive solutions if $\hat{p} < 0$ and $\hat{L}(X)$ is non-monotonic. Yet, the corresponding analysis can be applied in these cases alike.

5. *Estimation of Attractors Using Averaging Technique*

For systems with time dependent linear part, $p(t)$ frequently can be regarded as a highly oscillatory function, i.e., $p = p(t/\varepsilon)$, $0 < \varepsilon \ll 1$, where $\varepsilon$ is a characteristic time scale. Let us assume for simplicity that $k = k(t/\varepsilon)$, $L = L(t/\varepsilon, X)$, $F = F(t/\varepsilon)$, and $t_0 = 0$, and introduce fast time, $\tau = t/\varepsilon$, transforming (19) into the form,

(24)
$$D_\tau^+ X = \varepsilon \left( p(\tau) X + k(\tau) L(\tau, X) \right) + \varepsilon k(\tau) \|F(\tau)\| = \varepsilon q(\tau, X)$$
$$X(0) = X_0$$

The last equation assumes application of the averaging technique, which yields an autonomous counterpart to (24) that can be written as (20) under the following conditions:

(25)
$$\hat{p} = \varepsilon \lim_{T \to \infty} T^{-1} \int_0^T p(s) ds, \quad \hat{k} = \varepsilon \lim_{T \to \infty} T^{-1} \int_0^T k(s) ds, \quad \hat{L}(\|x\|) = \varepsilon \lim_{T \to \infty} T^{-1} \int_0^T L(s,\|x\|) ds,$$
$$\hat{F} = \varepsilon \lim_{T \to \infty} T^{-1} \int_0^T \|F(s)\| ds$$

Sufficient conditions assuring the closeness of some solutions of the averaged and original equations on large and infinite time-intervals can be found in [3], [14], [22], [23], and references therein. For instance, for $\tau \in (0,\infty)$ the following conditions assure closeness of solutions to (24) and (20) with coefficients defined by (25), see [14].

*Proposition*: Let the system (20) with coefficients defined by (25) have a positive solution, $X = d$, and $B_r(d)$ be a ball with radius $r$, which is centered at $X = d$. Assume that function $q(\tau, X)$ in (24) admits the following conditions:

1. $q(\tau, X) \leq \infty, \forall \tau \in \mathbb{R}, X \in B_r$, and the limit, $q_0(X) = \lim_{T \to \infty} T^{-1} \int_\tau^{\tau+T} q(s, X) ds$, exists uniformly $\forall X \in B_r$, and

   $\|\partial q(\tau, X)/\partial X\| \leq g_1 < \infty,$
   $\|dq_0(X)/dX\| \leq g_2 < \infty, \forall \tau \in \mathbb{R}, X \in B_r.$

2. $\partial^k q(\tau, X)/\partial X^k$ and $\partial^k q_0(X)/\partial X^k, k = 1, 2, \forall \tau \in \mathbb{R}, X \in B_r$ are continuous functions.

3. $q_0^k(X) = \lim_{T \to \infty} T^{-1} \int_\tau^{\tau+T} \left( \partial^k q(s, X)/\partial X^k \right) ds, k = 1,2$ are uniformly defined for $\forall \tau \in \mathbb{R}, X \in B_r$.

4. $dQ/dz(d) \neq 0$.

Then, for sufficiently small $\mu$, there is $\varepsilon_0 > 0$, such that for $0 < \varepsilon < \varepsilon_0$, (24) admits a unique solution $X^* = X^*(\tau)$, which obeys the inequality,

(26)
$$\sup \|X^*(\tau) - d\| < \mu, \forall \tau > 0$$

and $X^*(\tau)$ is the stable/unstable solution to (24) if $X = d$ is the stable/unstable solution of (20) with coefficients defined by (25).

Yet, positive $X_i^*(\tau)$-solutions belong to the $\mu_i$- neighborhood of the fixed solutions, $X_i = d_i$, to (20) with coefficients defined by (25) and adhere its stability/instability properties. Thus, $X_i^*(\tau)$ can be used for estimation of the trapping/stability regions of (24) the same way it was done prior with utility of the fixed solutions to (20) with coefficients defined by (21). Hence, for this case, (23) can be adjusted as follows: $\|W^{-1}(0)x_0\| \leq d_i - \mu$; hence the ellipsoids, $\omega_\gamma(0, d_i - \mu)|$, estimate the sets of pertained initial vectors for solutions to (1) or (2). This lets us adjust the statements of theorems 3 and 4 as follows:

**Theorem 5.** Assuming that (18) can be written in the form (24) with $\varepsilon \ll 1$, let functions $q(\tau, X)$ and $Q(X, \hat{F})$ admit the conditions of the Proposition, and:

1. $F = 0$, $\hat{p} < 0$, and (20) with coefficients defined by (25) possess a unique unstable fixed solution, $X = d$. Then, $x = 0$ solution to (2) is asymptotically stable. In addition, $\|x(t, x_0)\| \leq d + \mu$, $\forall x_0 \in \omega(0, d - \mu)$, $\forall t \geq 0$, and $\lim_{t \to \infty} \|x(t, x_0)\| \to 0$, $\forall x_0 \in \omega(0, d - \mu)$. If, in turn, $X = d$ is a stable unique fixed solution to (20) with coefficients defined by (25), then $\|x(t, x_0)\| \leq d + \mu$, $\forall x_0 \in \omega(t_0, d - \mu)$ and $\limsup_{t \to \infty} \|x(t, x_0)\| \leq d + \mu, \forall x_0 \geq 0$.

2. $\hat{F} \neq 0$, $\hat{p} < 0$, and (22) has two solutions, $d_1 > d_2$, then
$$\|x(t, x_0)\| \leq \begin{cases} d_1 + \mu, & \forall x_0 \in \omega(0, d_1 - \mu) \\ d_2 + \mu, & \forall x_0 \in \omega(0, d_2 - \mu) \end{cases}, \forall t \geq 0$$
$$\limsup_{t \to \infty} \|x(t, x_0)\| \leq d_2 + \mu, \forall x_0 \in \omega(0, d_1 - \mu)$$

3. $\hat{F} \neq 0$, $\hat{p} < 0$, (22) has a repeated root, $X = d$, and $x_0 \in \omega(0, d - \mu)$. Then, $\|x(t, x_0)\| \leq d + \mu$, $\forall x_0 \in \omega(0, d - \mu)$, $\forall t > 0$.

We notice that theoretical estimates for admissible values of $\mu$ and $\varepsilon(\mu)$ turn out to be quite conservative [14], but more accurate estimates frequently can be obtained in numerical simulations.

*Remark 2:* Application of averaging to (18), which includes two significantly different time scales, yields the equation possessing only slow-time, see [22] and [23] and more references therein. It was shown in [13] that under some conditions stability of the system averaged over fast-time implies stability of the original system with two-time scales. These inferences can be directly applied to (18) in the pertained cases. Moreover, after averaging over fast-time, slow-varying coefficients in (18) frequently can be effectively bounded, which lets us subsequently bound (18) by its time-invariant and integrable form.

*Remark 3:* Equation (18) turns into the integrable Bernoulli equation [24] if $F = 0$ and $L(t, X)$ obeys Holder's inequality, e.g. $L(t, X) \leq c(t) X^\alpha, c, \alpha > 0, X \in \Omega$ which streamlines stability analysis and estimation of the pertained solution bounds.

## 5. Simulations

This section applies the developed above methodology for estimating the solution norms as well as trapping/stability regions for Van der-Pol- like models with both time-varying linear part and external time-dependent perturbation. The system is written in dimensionless variables as follows,

(27)
$$\dot{x} = A(t)x + f(t, x) + F(t)$$
$$x(0) = x_0$$

where $x = [x_1\ x_2]^T$, $A = \begin{pmatrix} 0 & 1 \\ -\omega^2 & -\alpha_1 \end{pmatrix}$, $f = \begin{pmatrix} 0 \\ -\alpha_2 x_2^3 \end{pmatrix}$, and $F(t) = [0\ F_2(t)]^T$, $\omega^2(t) = \omega_0^2 + \omega_1(t), \omega_0 = const$, $\omega_1 = a_1 \sin r_1 t + a_2 \sin r_2 t$, $F_2 = a \sin \omega_2 t$, $a, \omega_2, a_i, r_i = const, i = 1, 2$. In all further simulations we set $\omega_0 = 2$ and $\alpha_1 = 0.2$.

Firstly, we notice that time-histories of $p(t)$, which are ubiquitous in our analysis, are affected by the normalization of $W(t)$. In particular, we contrast two different normalizations: 1). $W(0) = I$ and 2) $W(0) = W_0(0)$, where $W_0(t)$ is the fundamental matrix of solutions of the system $\dot{x} = Ax$ with $A = const$, i.e., $\omega_1 = 0$, and $\|W_0(0)\| = 1$.

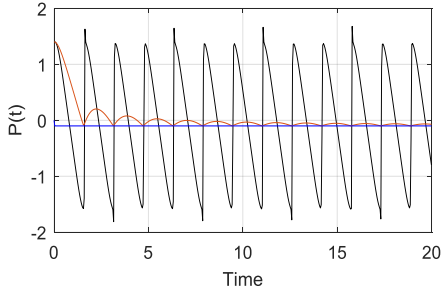

Fig.1. Black and blue lines plot time-histories of $p(t)$ computed for (27) in the following two cases: 1. $W(0) = I$, 2. $W(0) = W_0(0)$ respectively. Red line plots running average of $p(t)$ computed for the first normalization of $W(t)$.

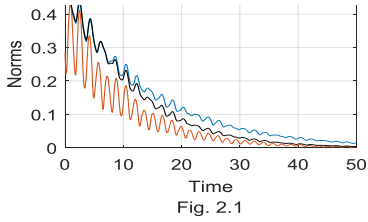

Fig. 2.1

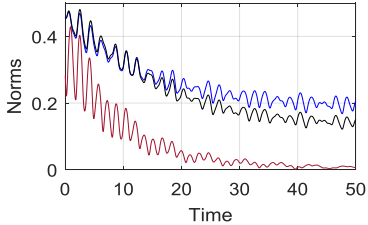

Fig. 2.2

Figure 2. Estimation of solution norms for system (27). Red, blue and black lines plot time-histories of the actual norm of solution to (27) and its upper bounds comprised by utility of Lipschitz or nonlinear Lipschitz inequalities, respectively.

In Fig. 1, black and blue lines plot time-histories of $p(t)$ corresponding to the first and second normalizations of $W(t)$, respectively. The red line plots running time-average of $p(t)$ for the first normalization of $W(t)$, i.e., $\bar{p}(t) = t^{-1} \int_0^t p(s)\,ds$. It is follow from the plots that the first normalization yields widely oscillatory $p(t)$, but the second – yields,

$p(t) = \mathrm{Re}(eig(A)) = -\alpha_1/2$, whereas $\lim_{t \to \infty} \bar{p}(t) = -\alpha_1/2$.

Further simulations show that the second normalization also reduces the variability of $p(t)$ if $|\omega_1| \leq \hat{\omega}$ and $\hat{\omega}$ assumes relatively small or intermediate values. Hence, the second normalization is adopted in our simulations.

The estimations of the norms of solutions to (27) are shown in Fig.2, where red, blue, and black lines plot time-histories of the actual norms of solutions and their two upper bounds, which are comprised by utility of linear (3) and nonlinear Lipschitz, i.e. $\alpha_2 x_2^3 \leq |\alpha_2| \|x\|^3$ inequalities, respectively. The value of Lipschitz constant depends upon $\sup x_2$ attaining in these simulations. This value is estimated using energy integral for the linearized, time invariant, and homogeneous model of (27). In these simulations, $\alpha_2 = 0.1,\ \omega_2 = 2\pi,\ a_1 = a_2 = 0.5,\ r_1 = \pi,\ r_2 = 7$; and in Fig. 2.1, $a = 0$, and in Fig. 2.2, $a = 0.01$.

Clearly, time – histories of the solution bounds comprising the nonlinear Lipschitz inequality outperform the ones utilizing (3) everywhere except a small initial time interval, where the latter is somewhat more accurate than former. Both bounds provide superior accuracy on the initial time intervals, which, however, decreases when time elapses. Bounds, based on (12), deliver tolerable accuracy on extended time intervals for the homogeneous system. However, the estimation accuracy declines for the nonhomogeneous system. We notice that the task of finding a suitable Lipschitz constant turns out to be rather deceptive for systems in higher dimensions. In contrast, devising a global nonlinear Lipschitz inequality is effortless for polynomial vector fields.

Fig. 3 contrasts the boundaries of trapping/stability regions, which are computed for (27) using three approaches: 1. simulation of (27), 2. simulation of the corresponding equation (18), and 3. utilization of approximate analytical models (20) with coefficients defined by either (21) or (25). Figs. 3.1-3.3 plot projections of the simulated solutions to (27) on $x_1 - x_2$ - plane, which approximate these three boundaries. Each approximation is defined by the initial vector with $x_2(0) = 0$, whereas for the first case, $x_1(0)$ is defined in simulations of (27). The pertaining solutions are shown in black lines. In the other two cases, $x_1(0)$ is defined in two successive steps. For both these cases, the first step defines the splitting values of $x_1(0)$. For approach 2, $x_1(0)$ is defined in simulations of the corresponding equation (18), whereas for approach 3, $x_1(0)$ is defined as the largest positive root of (20) with coefficients defined by either (21) or (25).

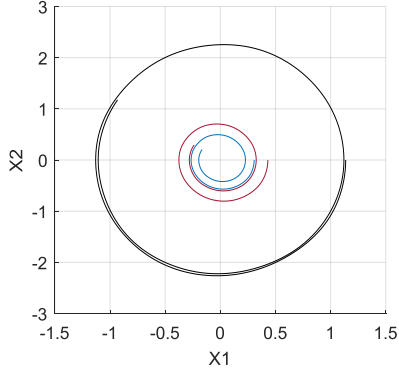

Fig. 3.1

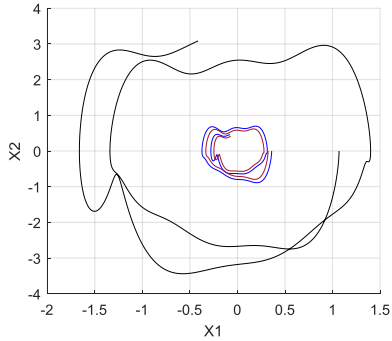

Fig. 3.2

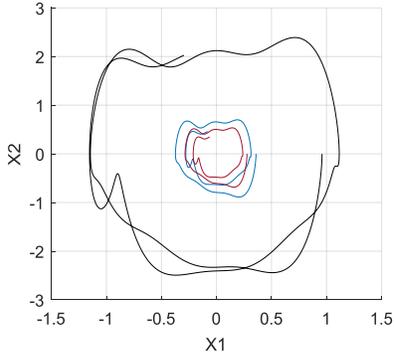

Fig. 3.3

Fig. 3. Estimation of stability basins and trapping regions. Figs. 3.1-3.3 plot projections of the solutions to (27) on $x_1 - x_2$ plane, which approximate the boundaries of the actual stability/trapping regions (black lines) and display two estimates to these boundaries in red and blue lines.

Subsequently, these values of $X_0$ are mapped into the values of $x_1(0)$ through application of the formula $\left\| W^{-1}(0) \begin{vmatrix} x_1(0) \\ 0 \end{vmatrix} \right\| = X_0$. The solutions to (27), which are emanated from these last two sets of the initial vectors, are plotted in red and blue lines, respectively.

All simulations in Fig. 3 and Fig. 4 are made with $\alpha_2 = -0.05$. The plots in Fig. 3.1 correspond to the following values of parameters, $a_1 = a_2 = 0.1$, $r_1 = 3.2\pi$, $r_2 = 13$, $a = 0.01$, whereas the plots in Figs. 3.2, 3.3 and Fig.4 are simulated for $a_1 = a_2 = 5$, $r_1 = 3.2\pi$, $r_2 = 13$, $\omega = 8\pi$. Note also that in Fig. 3.2, $a = 0$ and $\mu = 0.085$, and in Fig. 3.3 and Fig. 4, $a = 0.05$ and $\mu = 0.183$.

Clearly, simulations of (18) yield a central part of the actual stability basins or trapping regions for (27). We notice that the analytical approximations, which define the pertaining initial values through utility of (20) with coefficients defined by (21), are initially close to ones that are obtained in direct simulations of (18) if $\sup \omega_1(t)$ is relatively small and $a = 0$. If $a \neq 0$, the analytical approximations remain intact if $|a| \ll 1$, since in this case positive roots of (20) with coefficients defined by (21) sensitively depend upon $a$ and vanish if $a$ increases. Fig. 3.1 shows that under these last conditions the pertained values of $x_1(0)$, which are defined either by (20) with coefficients defined by (21) (blue lines) or in simulations of (18) (red lines), are close to each other. The difference between these two approximations to $x_1(0)$ further decreases if either $a = 0$ or $\omega_1 = 0$.

In turn, the averaging technique leads to tolerable analytical estimates of the trapping/stability regions for larger values of both $r_i$ and $a_i$. Fig. 3.2 and Fig. 3.3 compare numerical and analytical estimates for the pertained initial values that are developed in simulations of (18) and (20) with coefficients defined by (25). Clearly, the former two estimates, i.e., red and blue lines, are sufficiently close to each other for $a_i = 5$ and determine the central part of the actual trapping/stability regions of (27).

Fig.4 plots time-histories of $p(t)$ and $k(t)$, and their running time-averages in blue, yellow, red, and magenta lines, respectively. Both functions are notably oscillating, but their running time- averages quickly approach the constant values, which yield the principal contribution to the solutions of (18). Hence, these values can be used to derive the first approximation to the splitting values of $X_0$, which subsequently were refined in simulations of (18).

## 6. Conclusions and Future Work

This paper presents a novel approach to the analysis of stability, estimation of solution bounds, and trapping/stability regions of time-dependent and nonlinear systems, which is developed in the context of the first approximation methodology. This approach comprises the development of the pivoting differential inequality for the norm of solutions of the initial systems and subsequent analysis of the associated first order auxiliary differential equation. The linear and nonlinear parts of this equation are scaled by the functions, which are ubiquitous in our analysis. The time average of the former function approaches maximal Lyapunov exponent of the linear part of the initial system if $t \to \infty$ and, hence, determines its stability. The latter function turns out to be the running condition number of the fundamental matrix of the linearized system, which naturally scales its degree of robustness that varies in time.

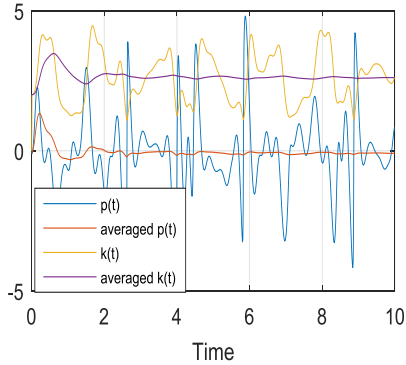

Fig. 4 plot time-histories of $p(t)$ and $k(t)$ in blue and yellow lines and their running time-averages in red and magenta lines, respectively.

We cast the auxiliary equation in the standard form by using either the Lipschitz inequality or its nonlinear extension; the latter is also charted in this paper. The application of Lipschitz condition linearizes the auxiliary equation and yields the corresponding solution bounds and stability criteria. In turn, adoption of the nonlinear extension of Lipschitz inequality leads to a more adequate nonlinear auxiliary equation, which extends our analysis beyond the standard framework of the first approximation methodology. We formulate the characteristic properties, which simplify numerical estimation of the trapping/stability regions of the nonlinear auxiliary equation and estimate the corresponding regions for the initial systems. Next, we introduce two approximations reducing the auxiliary equation to its autonomous and integrable forms. Consequently, we derive explicit estimates of stability/trapping regions for the auxiliary and initial systems and contrast the analytic and simulated results.

Our theoretical inferences are validated in inclusive numerical simulations that are partly presented in this paper. The simulations show that the accuracy of our estimates inversely correlates with the magnitudes of $\|f(t,x)\|$ and $\|F(t)\|$, since the auxiliary equation includes only the norms of these perturbations. Hence, the precision of the developed estimates turns out to be adequate if $\|f\|$ and $\|F\|$ are only known – a frequent premise in control theory. But, our approach can yield rather conservative estimates if both $f$ and $F$ are defined precisely.

Yet, the developed approach can be combined with some successive approximations to yield bilateral bounds for the norms of solutions that approach $\|x(t,x_0)\|$ under some broad conditions. Application of such refined methodology will be the topic of our subsequent paper.

This paper presumes that the fundamental solution matrix of the linearized system is known – a standard thesis in the first approximation methodology. This let to define key entries in our auxiliary equation, i.e. $p(t)$ and $k(t)$. For many practically important systems behavior of these functions on long time intervals can be evaluated through applications of pertained numerical simulations. Utility of analytical approximations to this matrix can strengthen our inferences and will be included in our future studies.

Acknowledgment. This paper was developed in close collaboration of its co-authors. The second co-author developed the programs and contributed to interpretation of the simulation results, whereas the first co-author developed the underlined methodology and primarily drafted the paper.